\documentclass[a4paper,10pt]{article}
\usepackage{amssymb}
\usepackage[T2A]{fontenc}
\usepackage[cp1251]{inputenc}
\usepackage{amsmath}

\usepackage[dvips]{graphicx}
\usepackage[english,russian]{babel}

\newtheorem{theorem}{Теорема}[section]

\newtheorem{definitionhead}[theorem]{Определение}

\author{А.~Я.~Канель-Белов}
\title{Памяти Валерия Анатольевича Сендерова}

\begin{document}
\maketitle

Я знал Валерия Анатольевича будучи школьником, он постоянно
организовывал математические бои, а я был их постоянным участником.
Помнятся несколько матбоев, где Саша Разборов был капитаном команды
Второй школы, я заместителем, а капитаном команды 91-й школы был
Максим Концевич.

Сендеров
научил меня
одной важной
вещи. Говоря о решении задач, он показывал идейное
ядро, где всё и происходит. Оно маленькое -- это как жало станка, и
именно это --
главное, что надо увидеть. Станок состоит из большой станины, приводных
ремней
таблички с указанием производителя,
и т.п.,
а жало маленькое. Так же и
задача: важно выделять, где всё происходит и почему происходит.

Выделению ядра или ``жала'' он научил не только меня. Привожу воспоминание одного
математика:

\medskip

{\it\small Мне тоже вспомнилось что-то из моего небольшого опыта общения с
В.А. Так случилось, что вышел он из тюрьмы в  1987 году, когда я
был в выпускном классе.
 И возникла
активность по натаскиванию всех желающих
из Второй и 57-й школ на мехматские ``гробы'', с тем чтобы поступить на мехмат,
точнее поступать.
Я на мехмат поступать не стал, в чем ни на минуту не
раскаиваюсь,
 но вот сами эти занятия \dots Их было несколько. Разбирались
реально
сложные задачи, вполне олимпиадные, но в
олимпиадных предполагается элемент красоты, что для мехматских
``гробов'' необязательно.
И вот, хоть я никогда олимпиады особо не любил, эти занятия
имели какую-то особую красоту.
Действительно,
как Лёша пишет,
В.А. показывал некоторое ``ядро'' и ``всё остальное'', и делал это
так мастерски, как никто другой в контексте олимпиадных задач. Я первый
раз в жизни ощутил какую-то красоту олимпиадных задач, и мне хотелось
ходить на эти занятия независимо от мехмата и т.п. И это была даже
красота не самих задач, а красота раскладывания на ``ядро'' и ``всё
остальное''. Ну и, конечно, было ощущение, что ты общаешься с совершенно
героическим человеком, просидевшим 5 лет в тюрьмах, значительную часть
этого в карцере. В.А. несомненно обладал некоторым гипнотизмом,
исходившим из его могучей внутренней силы. Людям, окружавшим меня,
от родителей до Гриши К.,
потребовалось немало усилий, чтобы заставить меня
рассмотреть незамутнённым сознанием слова В.А. ``ну если вас завалят,
что наверняка и будет, то по вашим костям пройдут другие''.
Когда я находился рядом с ним,  критическое мышление
(типа ``а завалюсь так заберут в армию'' и проч.) совсем отшибало.
Этот человек пошел на намного бОльшие
лишения, и думать о собственном комфорте рядом с ним было невозможно.
Я слышал от разных людей, особенно старшего поколения, с которыми В.А.
общался до тюрьмы, что так реагировали очень многие.
}\footnote{ Автор не одобряет вовлечение молодежи, особенно несовершеннолетних, в политику}

\medskip

Валерий Анатольевич Сендеров вместе с Борисом Ильичом Каневским заложили традицию олимпиад и математических боев во Второй школе, продолжающуюся и по сей день.
 Очень часто эти матбои Вторая Школа выигрывала. Такая
традиция принесла плоды не только во Второй школе. В последующем покойный Митя Дерягин (выпускник 1981 года, победитель Всесоюзной олимпиады) начал кодификацию правил. Более-менее окончательную форму матбои приняли в начале 90-х, после синтеза московской и питерской версии правил, большая в этом  заслуга Саши Ковальджи (ныне -- зам. директора по науке лицея ``Вторая Школа''). \footnote{См. 
Дерягин Д.~В., Канель А.~Я., Ковальджи А.~К., Кондаков Г.~В., Рубанов И.~С., Финашин С.~М., Фомин Д.~В., Шапиро А.~А., Яценко А.~Д, ``Математический бой двух команд: Правила, комментарии, опыт проведения'', Математика в школе, 1990, № 4, 20--25.}
Упомянем Московские Турниры математических боев.

Чтобы не создалось искажённого представления об общественных взглядах Валерия Анатольевича, следует отметить, что он был государственник. Его последнее интервью можно найти по ссылке http://www.russ.ru/pole/Kak-byvshij-dissident-i-politzaklyuchionnyj-stanovitsya-ohranitelem . Заголовок и манера, в которой было взято интервью, автору не нравятся -- в конце концов, интервьюер не должен давать ярлыки, тем более в заголовке. Автор приводит эту ссылку только потому, что это интервью -- последнее. Упомянем также последнюю статью Валерия Анатольевича
(совместно с Ю.~Кублановским и Ф.~Разумовским) http://www.rg.ru/2014/03/12/pismo.html.
Он высоко ценил  ``Вехи''~-- сборник статей о русской интеллигенции, созданный деятелями Серебрянного Века. Этот сборник перевернул моё сознание,  дал мне порцию свободы и понимания, в том числе того слоя людей, с которым часто приходится иметь дело. Одна из мыслей ему и мне, как человеку занимающегося классификацией идей решения олимпиадных задач, оказалась близка. Давайте разберемся, что такое ``правые'' и что такое ``левые''. Для этого выпишем типичные ``правые'' и типичные ``левые'' взгляды и заметим, что люди мыслят ``пакетно''. (Вспоминается остроумное высказывание Валерия Анатольевича о характере дискуссий: ``они каются в грехах друг друга'').  Эти взгляды надо объяснять, но не исходя из их ``истинности'' или ``ложности'', но исходя из эмоционального строя. И такая попытка, пусть весьма не полная, в сборнике была дана. Я бы добавил, что анализ должен использовать технику, в частности, З.~Фрейда и К.~Юнга, находивших ``смысл'' в симптомах. И в этом есть родство с книгой В.~Проппа ``Исторические корни волшебной сказки'' (в первой части классифицируются сюжеты и элементы сказки, дается структурный анализ, во второй даются объяснения).

Есть люди, с которыми не во всём соглашаешься, но с их уходом возникает некая пустота (такими, например, для
автора были Н.~Б.~Васильев и И.~Ф.~Шарыгин). Я хотел обсудить
с В.А. ряд вешей, собирался позвонить
-- но как-то всё откладывалось...

\end{document}